\documentclass[12pt]{article}
\usepackage{amsfonts}
\usepackage{mathptmx}
\usepackage[utf8]{inputenc}
\usepackage{amsmath}
\usepackage{breqn}
\usepackage{mathrsfs}
\newtheorem{thm}{Theorem}
\newtheorem{defn}{Definition}
\newtheorem{lem}{Lemma}
\def\ff{ \mathbf{F} }

\title{Higher Level Completeness for Permutation Polynomials}
\author{S Rajagopal  and P Vanchinathan*}
\begin{document}
\maketitle
\begin{center}
Division of Mathematics\\
VIT University\\
Vandalur–Kelambakkam Road\\
Chennai, 600 127 INDIA\\[3pt]
\texttt{rajagopal.s2020@vitstudent.ac.in,\quad 
		vanchinathan.p@vit.ac.in}
\begin{abstract}
Generalising the concept of a complete permutation polynomial over 
a finite field, we define
completness to level $k$ for $k\ge1$ in fields of odd characteristic.
We construct  two families of polynomials that satisfy the condition 
of high level completeness for all finite fields, and
two more families complete to the  maximum level possible for 
large collection of finite fields.

Under the binary operation of composition of functions
one family of polynomials is an abelian group
isomorphic to the additive group, while the other is 
isomorphic to the multiplicative group.
\end{abstract}

\end{center}

\section{Introduction}
A polynomial over a finite field is said to be a permutation polynomial
if it gives rise to a bijective function on that field.
The monograph \cite{lidl} has a whole chapter devoted to this.
The handbook has a large collection of results and bibliography
\cite{panario}. The survey articles by Hou \cite{hou1}, \cite{hou2}, gives
the status of the theory of permutation in general and 
binomials and trinomials in particular.
Due to its applications to combinatorics, cryptology and coding theory
lot of effort is put in constructing them.

 A  permutation polynomial $f(x)\in\ff_q[x]$ is  called a
complete permutation polynomial, if $f(x)+x$ is also
a permutation polynomial.

Complete permutation polynomials seem to have certain
advantages.   In applications to cryptography
they are used for constructing bent functions.
Results on complete permutation polynomials can be found,
 for example, in \cite{bart1},\cite{bart2}, \cite{bass}, \cite{cao}, \cite{tu}.

In this paper our objective is to find large family of complete
permutation polynomials in each finite field and we achieve it 
and as a bonus, they are higher level complete polynomials:

\begin{defn} For a positive integer $k$, we say a polynomial 
$f(x)\in \ff_q[x]$ is complete to level $k$, or simply $k$-complete, 
if along with $f(x)$ we  have $f(x)+x, f(x)+2x, \ldots, f(x)+kx$ 
 are also permutation polynomials.
\end{defn}
Clearly complete polynomials of level upto $p-1$
where $p$ is the characteristic of $\ff_q$, only need to be studied. 

Here are some easy-to-verify statements of higher level
completeness property of a permutation polynomial:
\begin{itemize}
\item Higher completeness is same as usual completeness in
characterstic $2$.
\item Usual complete polynomials are 1-complete.
\item If a permutation polynomial is $k$-complete,
then it is also $k'$-complete for positive integers $k'<k$.
\end{itemize}
\begin{defn}
A permutation polynomial that is complete to level $k$  for $k=p-1$
is said to be a  maximally complete permutation polynomial.
\end{defn}

\noindent\textbf{Example:} 
Any linear polynomial $ax+b\in \ff_q[x]$ is maximally complete
 if $a$ is not in the prime subfield; in case $a$ is in the 
 prime subfield, it will be complete to level  $k=p-1-a' $,   
with $a'$ being the least positive integer representing $a\in \ff_p$.

\vspace{6pt}
\noindent \textbf{Question}:\quad \textit{Do complete permutation
polynomials of level $k>1$ and degree${}>1$ exist?}

\vspace{6pt}
Our main results of this paper are the following
affirmative  answers to the above question:

\vspace{6pt}
\noindent\textbf{Theorem A}:\quad
\textit{For any prime $p\ge3$ and a composite integer $n$, finite
fields of order $p^n$ admit maximally complete permutation polynomials
	which are non-linear.}

\vspace{6pt}
In case the finite field does not have a proper subfield bigger than the
prime subfield (the case where $n$ in Theorem A is a prime)
we can find polynomials that are almost maximally complete:

\vspace{6pt}
\noindent\textbf{Theorem B}:\quad
\textit{In  all finite fields of characteristic $\ge3$, there 
exist non-linear complete permutation polynomials of level $p-2$.}

\vspace{6pt}
\textit{For a more precise and constructive version of the above, 
see Theorems~\ref{Completeppadd}, \ref{Completeppmul} and
\ref{maxcpp}} in the next section.

\section{Main Results and Proofs}
Notation:\quad  $ \ff_{q}$ will denote a finite field of  
prime power order $q$, assumed to be of odd characteristic $p$.
 
For an extension field  $\ff_{q^n}$  of $ \ff_{q}$.
we denote by   $ m = q + q^2 + q^3  +\cdots+ q^{n-1}$. 
Note that $m+1 = q^n-1$ and that $m$ is  a multiple of 
the characteristic $p$. 

\begin{lem}  \label{addpp}
For a finite field $\ff_q$ and an extension field $\ff_{q^n}$
and any $ c \in \ff_{q}  $ the following polynomial
 $$ f_{c+} (x)  = x+c \sum \limits _{j=1}^{m} x^{j(q-1)}$$ 
is a permutaion polynomial over $\ff_{q^n}$.
\end{lem}

Proof: It will follow once we  show that the  polynomial 
$f_{c+}(x)$ actually is the  function given below:
\begin{equation*}
f_{c+}(a)= \left\{ \begin{array}{lll}  a -c &  \mbox{for } & a \in \ff_{q^n}\setminus \ff_{q} \\ 
 a & \mbox{for } & a \in \ff_{q}
\end{array}\right.
\end{equation*}

When $a$ is in the base field $\ff_q$, we claim that all the  
terms in the summation are 1. That is so  because we are 
summing the $r$-th powers of $a$ where $r$ is always
a multiple of $q-1$. Thus $f_{c+}(a) = a + cm = a$.

Next we will move on to evaluating it on the  remaining elements. 
So consider $a\in \ff_{q^n}\setminus \ff_q$. The summation  
actually is a sum of $m$ terms of a geometric progression. So
\begin{align*}
f_{c+}(a)  & =   a + c \left[a^{(q-1)}  +  a^{2(q-1)}  + \cdots + a ^{m(q-1)} \right]\\
       &= a  + c a^{q-1}\Biggl[ \frac{\left( a^{q-1} \right) ^ m -1 }{a^{q-1}-1}
         \Biggr]\\
        &= a  + c\Biggl[\frac{a^{q^n-1}-a^{q-1}}{a^{q-1}-1}\Biggr]
	\mbox{ as } (m+1)(q-1) = q^n-1\\
        &= a  + c\Biggl[\frac{1-a^{q-1}}{a^{q-1}-1} \Biggr] \\
	& = a  - c      \hspace{2cm} \square
\end{align*}

Summarising we see that $f_{c+}(x)$ when restricted to the base 
field  is the identity function, and, on the elements not in 
the base field, is  a translation by an element of the base field.
Thus $f_{c+}(x)$ is a permutation polynomial of $\ff_{q^n}$ for
all $c\in \ff_q$.
\textit{Though we make no use of it further, it is now
 clear that the permutation given by this polynomial for $c\ne0$, 
 has order $p$, with exactly $q$ fixed points}.

Next we provide a multiplicative analogue of the polynomials 
constructed in the previous theorem.
\begin{lem}\label{multpp}
Under the same notation as earlier,  we define
for $ c \in \ff_{q} , c\neq 1 $ a polynomial by	
$$ f_{c*} (x)  =  x +cx \sum \limits _{j=1}^{m} x^{j(q-1)}. $$ 
This  is a permutaion polynomial over $\ff_{q^n}$.
\end{lem}
In fact this polynomial is the function given by the following description:
\begin{equation*}
f_{c*}(x)= \left\{ \begin{array}{lll}  x \left(1-c\right) &  \mbox{for } & x \in \ff_{q^n}\setminus \ff_{q} \\
 x & \mbox{for } & x \in \ff_{q}
\end{array}\right.
\end{equation*}

\noindent
\textbf{Case 1 :} \\
      To prove  $ f_{c*}(a) =  a   \forall a  \in  \ff_{q}  $. \\
       As $f_{c*}(0)=0$, we consider $a\in \ff_q^*$.
\begin{align*}
&f_{c*}(a) = a + c\left[a^{q} + a^{(2q-1)} +\cdots+ a^{m.q-(m-1)} \right] \\
& =  a + ca\big[a^{q-1} + a^{2(q-1)}+ a^{3(q-1)}+ \cdots + a^{m(q-1)}  \big]\\
& =  a + ca\left[\underbrace{a +a+ \cdots + a}_{m\rm\ times} \right]\\
& =  a + cma = a
\end{align*}

\noindent \textbf{Case 2 :} $a\in \ff_{q^n}\setminus \ff_{q}$.
\begin{align*}
f_{c*}(a) 
	& = a + ca \left[a^{(q-1)} + a^{2(q-1)}+ a^{3(q-1)}
	        + \cdots + a^{m(q-1)}  \right]\\
	& = a + ca\Biggl[\frac{a^{(m+1)(q-1)}-a^{q-1}}{a^{q-1}-1}\Biggr]\\
        & = a + ca\Biggl[\frac{1-a^{q-1}}{a^{q-1}-1} \Biggr] \\
        & = a\left(1-c\right)   \hspace{2cm}\square
\end{align*}

\textit{Again we can specify the cycle type of this
permutation $f_{c*}$ in terms of the order of $(1-c)$ as an element
of the multiplicative group} $\ff_q^*$.

The permutation polynomials constructed above are actually
complete to higher level. Here is the precise result:
\begin{thm}\label{Completeppadd}
The polynomial $f_{c+}(x)$ is $(p-2)$-complete,
that is, $f_{c+}(x)$ along with the set of $(p-2)$ 
polynomials $x+f_{c+}(x),2x+f_{c+}(x),3x+f_{c+}(x),\ldots,
	(p-2)x+f_{c+}(x)$ are all  permutation polynomials.
\end{thm}
Proof: For a positive integer $k$ from the definition we get 
\begin{align*}
kx + f_{c+}(x) 
	& = (k+1)x + c\left[x^{q-1} + x^{2(q-1)}
        	   + \cdots + x ^{m(q-1)} \right]\\
	& = (k+1) \Biggl[x + \frac{c}{k+1} \left(x^{q-1} 
 	         +  x^{2(q-1)} +\cdots + x ^{m(q-1)} \right)\Biggr]
\end{align*}
But, visual inspection makes it clear that,
the expression inside the square brackets is actually
$f_{c'+}$ for $c'= c/(k+1)$. 
Thus we arrive at an identity 
\begin{align*}
   kx+ f_{c+}(x) & = \left(k+1\right)f_{c'+} (x)
\end{align*}
By Lemma~\ref{addpp} RHS is a permutation polynomial and hence
because the LHS is too. Of course we need  $k+1\ne0$, and so
this is true for $k=1,2,\ldots,p-2$.
So $f_{c+}(x)$ is $(p-2)$-complete. \hfill $\square$
\begin{thm}\label{Completeppmul}
The polynomial $f_{c*}(x)$ are complete to level $p-2$, i.e., 
$f_{c*}(x)$ along with the set of $(p-2)$ polynomials
$x+f_{c*}(x),2x+f_{c*}(x),3x+f_{c*}(x),\ldots,(p-2)x+f_{c*}(x)$
are  permutation polynomials.
\end{thm}
Proof: Using arguments similar to the ones employed in the proof
of Theorem~\ref{Completeppadd} we can  easily arrive at
\begin{align*}
kx + f_{c*}(x) & = \left(k+1\right)f_{c'* } (x)
\end{align*}
for $c'= c/(k+1)$.
Again we will be able to conclude  that
$kx + f_{c*}\left(x\right) $  is also permutation 
polynomial for $ 0< k<(p-1)    $.  \hfill $\square $

Next we move on to finding maximally complete permutation polynomials.
First we need a temporary definition of a \textit{middle  subfield}
of  a field:  by that we mean a proper subfield which is not the prime subfield.
Middle subfields exist in finite fields of order $p^n$ iff $n$ is a composite
number.
\begin{thm}\label{maxcpp}
Let $\ff_q$ be a finite field admitting at least one middle  
subfield. Choose $b,c\in\ff_q$ such that $b$ is in some middle 
subfield and  $c$ is in any proper subfield of $\ff_q$
Then the polynomials $bf_{c+}(x)$ and $bf_{c*}(x)$ are
 maximally complete permutation polynomials.
\end{thm}

Proof: 
As  $c$ is  in a proper subfield of $\ff_q$,  Lemma~\ref{addpp}
ensures $f_{c+}(x)$ is a permutation polynomial
and hence the scalar multiple $bf_{c+}(x)$ is  too.
We will now show for any $k=1,2,\ldots, p-1$ 
that $bf_{c+}(x) +kx$ is a permutation polynomial.

We simply rewrite this as 
$$bf_{c+}(x)+kx =  (b+k) \bigg(x +\frac{bc}{b+k}\sum_j x^{j(q-1)}\bigg)
$$
And the latter polynomial is the same as $(b+k)f_{c'+}(x)$ 
where $c'= bc/(b+k)$, and hence a permutation
polynomial. We need $b+k$ to be nonzero for this to be true.
That follows from the hypothesis that $b$ is not in the prime 
field while $k$ is always in the prime field. 
This completes the proof that $bf_{c+}(x)$ is maximally complete.

Same arguments work for proving maximal completeness
property for the polynomials $bf_{c*}$ and so omitted.

Next we discuss inter- and intra-relationship among the
members of two families of polynomials $f_{c+}(x)$ and $f_{c*}(x)$.

\begin{thm}\label{groupstructure}
The two collections of functions $f_{c+}(x)$, and 
$f_{c*}(x)$ behave well under composition. In fact, 
\begin{itemize}
\item[(i)]  $\{\, f_{c+}(x)\mid c\in \ff_q\, \}$ is an abelian group under
	composition, and is  isomorphic to the additive group of $\ff_q$.
\item[(ii)] $ \big\{\,f_{c*}(x)\mid c\in \ff_{q}-\{1\}\, \big\} $ is also an abelian group under composition, and is isomorphic to the multiplicative
	group $\ff_q^*$.
\item[(iii)] $ f_{c+}(x) $ and $ f_{c*}(x)$ have the following relationship: 
$$x  \big(f_{c+}(x) -x+1\big )  = f_{c*}(x)$$
	\end{itemize}
\end{thm}
	Proof:

Verifying (i) is straightforward.

For part (iii) the  given relationship  between $f_{c+}(x)$ 
and $f_{c*}(x)$ follows directly from their definitions.

For (ii), to see that it is a group one can easily check that
\begin{align*}
   f_{c*} \circ f_{d*} (x) &=f_{(c+d-cd)*}
	= \left\{ \begin{array}{lll}  
		(\,1-c-d+cd\,)x 
		&  \mbox{for } & x \in \ff_{q^n}\setminus \ff_{q} \\
              x & \mbox{for } & x \in \ff_{q}
\end{array}\right.
\end{align*}
 and that the compositional inverse of $f_{c*}(x)$ is 
\begin{align*}
f_{c*}^{-1}(x)  & = \left\{ \begin{array}{lll} c x/(c-1)
        	& \mbox{for } & x \in \ff_{q^n}\setminus \ff_{q} \\
              x & \mbox{for } & x \in \ff_{q}
\end{array}\right.
\end{align*} 
Now  the statement that is it isomophic to $\ff_q^*$
can be deduced from the following lemma whose proof
is a simple exercise:
\begin{lem}
For any field $K$,  the set $K\setminus\{1\}$
is a group under the binary operation
$a*b = a+b-ab$ and this is isomorphic to $K^*$
under the map $a\mapsto 1-a$.
\end{lem}

\vspace{1pc}
\textbf{Examples:}
First two polynomials given below are 3-complete permutation polynomials
over $\ff_{25}$, and the other are 5-complete over $\ff_{49}$.
\begin{itemize}
\item $ f_{2+}(x) = 2x^{20} + 2x^{16} + 2x^{12} + 2x^8 + 2x^4 + x \in \ff_{5^2}[x]$ \\
\item $f_{4*}(x) = 4x^{21} + 4x^{17} + 4x^{13} + 4x^{9} + 4x^{5} + x\in
	\ff_{5^2}[x]$ \\
\item $f_{4+}(x)= 4x^{42}+4x^{36}+4x^{30}+4x^{24}+4x^{18}+4x^{12}+4x^{6}+x \in \ff_{7^2}[x]$\\
\item $f_{6*}(x)= 6x^{43} + 6x^{37} +6x^{31} + 6x^{25} + 6x^{19} + 6x^{13} + 6x^{7} + x \in \ff_{7^2}[x]$ 
\end{itemize}

\end{document}